\baselineskip=17pt plus3pt minus1pt 
\magnification=\magstep1       
\hsize=5.5truein                     
\vsize=8.5truein                      
\parindent 0pt
\parskip=\smallskipamount
\mathsurround=1pt
\hoffset=.25truein
\voffset=.5truein
\nopagenumbers
%
%
\def\today{\ifcase\month\or
  January\or February\or March\or April\or May\or June\or
  July\or August\or September\or October\or November\or December\fi
  \space\number\day, \number\year}
%
 at 10truept

%
\newcount\dispno      
\dispno=1\relax       
\newcount\refno       
\refno=1\relax        
\newcount\citations   
\citations=0\relax    
\newcount\sectno      
\sectno=0\relax       
\newbox\boxscratch    
%

%
%
%
\def\Section#1#2{\global\advance\sectno by 1\relax%
\label{Section\noexpand~\the\sectno}{#2}%
\smallskip
\goodbreak
\setbox\boxscratch=\hbox{\bf Section \the\sectno.~}%
{\hangindent=\wd\boxscratch\hangafter=1
\noindent{\bf Section \the\sectno.~#1}\nobreak\smallskip\nobreak}}
%
\def\sqr#1#2{{\vcenter{\vbox{\hrule height.#2pt
              \hbox{\vrule width.#2pt height#1pt \kern#1pt
              \vrule width.#2pt}
              \hrule height.#2pt}}}}
\def\square{$\mathchoice\sqr34\sqr34\sqr{2.1}3\sqr{1.5}3$}
\def\endproof{~~\hfill\square\par\medbreak}
\def\noproof{~~\hfill\square}
%
%
\def\proc#1#2#3{{\hbox{${#3 \subseteq} \kern -#1cm _{#2 /}\hskip 0.05cm $}}}
\def\propcont{\mathchoice\proc{0.17}{\scriptscriptstyle}{}
                         \proc{0.17}{\scriptscriptstyle}{}
                         \proc{0.15}{\scriptscriptstyle}{\scriptstyle }
                         \proc{0.13}{\scriptscriptstyle}{\scriptscriptstyle}}
%

%
\def\normalin{\hbox{\raise0.045cm \hbox
                   {$\underline{\triangleleft }$}\hskip0.02cm}}
%
%
\def\'#1{\ifx#1i{\accent"13 \i}\else{\accent"13 #1}\fi}
%
%
%
\def\semidirect{\rlap{$\times$}\kern+7.2778pt \vrule height4.96333pt
width.5pt depth0pt\relax\;}
%
%
\def\prop#1#2{\noindent{\bf Proposition~\the\sectno.\the\dispno. }%
\label{Proposition\noexpand~\the\sectno.\the\dispno}{#1}\global\advance\dispno 
by 1{\it #2}\smallbreak}
\def\thm#1#2{\noindent{\bf Theorem~\the\sectno.\the\dispno. }%
\label{Theorem\noexpand~\the\sectno.\the\dispno}{#1}\global\advance\dispno
by 1{\it #2}\smallbreak}
\def\cor#1#2{\noindent{\bf Corollary~\the\sectno.\the\dispno. }%
\label{Corollary\noexpand~\the\sectno.\the\dispno}{#1}\global\advance\dispno by
1{\it #2}\smallbreak}
\def\defn{\noindent{\bf
Definition~\the\sectno.\the\dispno. }\global\advance\dispno by 1\relax}
\def\lemma#1#2{\noindent{\bf Lemma~\the\sectno.\the\dispno. }%
\label{Lemma\noexpand~\the\sectno.\the\dispno}{#1}\global\advance\dispno by
1{\it #2}\smallbreak}
\def\rmrk#1{\noindent{\bf Remark~\the\sectno.\the\dispno.}%
\label{Remark\noexpand~\the\sectno.\the\dispno}{#1}\global\advance\dispno
by 1\relax}
\def\proof{\noindent{\it Proof: }}
\def\numbeq#1{\the\sectno.\the\dispno\label{\the\sectno.\the\dispno}{#1}%
\global\advance\dispno by 1\relax}

\def\comm#1,#2{\left[#1{,}#2\right]}
\newdimen\boxitsep \boxitsep=0 true pt
\newdimen\boxith \boxith=.4 true pt 
\newdimen\boxitv \boxitv=.4 true pt
\gdef\boxit#1{\vbox{\hrule height\boxith
                    \hbox{\vrule width\boxitv\kern\boxitsep
                          \vbox{\kern\boxitsep#1\kern\boxitsep}%
                          \kern\boxitsep\vrule width\boxitv}
                    \hrule height\boxith}}
\def\square{\ \hbox{\vrule height7.5pt depth1.5pt width 6pt}\par}
\outer\def\square{\ifmmode\else\hfill\fi
   \setbox0=\hbox{} \wd0=6pt \ht0=7.5pt \dp0=1.5pt
   \raise-1.5pt\hbox{\boxit{\box0}\par}
}

\def\frac#1/#2{\leavevmode\kern.1em
              \raise.5ex\hbox{\the\scriptfont0 #1}\kern-.1em
              /\kern\.15em\lower.25ex\hbox{\the\scriptfont0 #2}}
\def\incnoteq{\lower.1ex \hbox{\rlap{\raise 1ex
     \hbox{$\scriptscriptstyle\subset$}}{$\scriptscriptstyle\not=$}}}
%
%


\def\propcontup{\bigcup\!\!\!\rlap{\kern+.2pt$\backslash$}\,\kern+1pt\vert}
%
%
%
\def\label#1#2{\immediate\write\aux%
{\noexpand\def\expandafter\noexpand\csname#2\endcsname{#1}}}
%
\def\ifundefined#1{\expandafter\ifx\csname#1\endcsname\relax}
%
%
\def\ref#1{%
\ifundefined{#1}\message{! No ref. to #1;}%
 \else\csname #1\endcsname\fi}
%
%
\def\refer#1{%
\the\refno\label{\the\refno}{#1}%
\global\advance\refno by 1\relax}
%
%
\def\cite#1{%
\expandafter\gdef\csname x#1\endcsname{1}%
\global\advance\citations by 1\relax
\ifundefined{#1}\message{! No ref. to #1;}%
\else\csname #1\endcsname\fi}
%
%
 at 8truept      
%
%
%

\newread\aux
\immediate\openin\aux=\jobname.aux
\ifeof\aux \message{! No file \jobname.aux;}
\else \input \jobname.aux \immediate\closein\aux \fi
\newwrite\aux
\immediate\openout\aux=\jobname.aux
 
 at 8truept

\null\vskip 1truein
\centerline{DOMINIONS IN FINITELY GENERATED NILPOTENT GROUPS}
\vskip \baselineskip
\vskip \baselineskip
\centerline{Arturo Magidin}
\vskip \baselineskip
{\baselineskip 12pt plus2pt minus1pt
\leftskip 2.2truein Oficina 112\par
\leftskip 2.2truein Instituto de Matem\'aticas\par
\leftskip 2.2truein Circuito Exterior\par
\leftskip 2.2truein Ciudad Universitaria\par
\leftskip 2.2truein 04510 M\'exico City, MEXICO\par
\leftskip 2.2truein e-mail: magidin@matem.unam.mx\par}
\vskip \baselineskip
\vskip \baselineskip
\vskip \baselineskip
\vskip \baselineskip

{\parindent=20pt
\narrower\narrower
\noindent{{Abstract. In the first part, we prove that
the dominion (in the sense of Isbell) of a subgroup of a finitely
generated nilpotent group is trivial in the category of all nilpotent
groups. In the second part, we show that the dominion of a subgroup of
a finitely generated nilpotent group of class~two is trivial in the
category of all metabelian nilpotent~groups.\par}}}
\bigskip
\medskip



\Section{Introduction}{intro}

Suppose that a nilpotent group~$G$ and a subgroup~$H$ of~$G$ are
given. Are there any elements $g\in G\setminus H$ such that any two
group morphisms to a nilpotent group which agree on~$H$ must also
agree on~$g$?

To put this question in a more general context, let~$\cal C$ be a full
subcategory of the category of all algebras (in the sense of Universal
Algebra) of a fixed type, which is closed under passing to
subalgebras. Let $A\in {\cal C}$, and let~$B$ be a subalgebra
of~$A$. Recall that, in this situation, Isbell {\bf
[\cite{isbellone}]} defines the {\it dominion of~$B$ in~$A$} (in the
category~${\cal C}$) to be the intersection of all equalizer
subalgebras of~$A$ containing~$B$. Explicitly,
$${\rm dom}_A^{\cal C}(B)=\Bigl\{a\in A\bigm| \forall C\in {\cal C}\;
\forall f,g\colon A\to C,\ {\rm if}\ f|_B=g|_B{\rm\ then\ }
f(a)=g(a)\Bigr\}.$$

Note that the dominion depends on the category of context, that it is
a subalgebra of~$A$, and that it
always contains~$B$. If $B\propcont {\rm dom}_A^{\cal C}(B)$, we say
that the dominion of~$B$ in~$A$ (in the category~${\cal C}$) is
{\it nontrivial\/}, and we say the dominion is {\it trivial\/}~otherwise.

Therefore, the question above is equivalent to asking whether the
dominion of~$H$ in~$G$ (in the category~${\cal N}\!il$, consisting of
all nilpotent groups) is~nontrivial.

To provide some context for the results in this paper, we note that
in~{\bf [\cite{nildomsprelim}]} we proved that there exists an
infinitely generated nilpotent group~$G$ of class~two, and a finitely
generated subgroup~$H$ of~$G$, such that
$$H\propcont {\rm dom}_G^{{\cal N}_2}(H)={\rm dom}_G^{{\cal
N}\!il}(H)$$
(where ${\cal N}_2$ denotes the variety of all nilpotent groups of
class at most~$2$, and ${\cal N}\!il$ the category of all nilpotent groups),
and the dominion of~$H$ in~$G$ is not finitely generated. We also
proved that there exists a {\it finitely} generated nilpotent
group~$G$ of class~two, such that for any fixed $c>1$ there exists a
subgroup $H_c$ of~$G$ such that
$$H_c\propcont {\rm dom}_G^{{\cal N}_2}(H_c)={\rm dom}_G^{{\cal N}_c}(H_c),$$
where ${\cal N}_c$ is the category of all nilpotent groups of class at
most~$c$. We also indicated there, but gave no proof, that dominions of
subgroups of finitely generated nilpotent groups are trivial in~${\cal
N}\!il$, filling the gap between the two results quoted above. We will
provide a proof of this fact in the present~work.

In \ref{notation} we will review basic facts about dominions which we
will use in later parts, as well as establish notation.

The main results of this work are in two parts; in the first part,
\ref{generalfg}, we will prove that if~$G$ is a finitely generated
nilpotent group (of any class), then for all subgroups~$H$ of~$G$ we
have that ${\rm dom}_G^{{\cal N}\!il}(H)=H$.
This section is elementary, and only assumes knowledge of basic facts about
finitely generated nilpotent groups, and a result of~Higman {\bf
[\cite{higmanpgroups}]} on amalgams of~$p$-groups. We will recall these
results~below.

In the second part, \ref{generalforniltwo}, we will prove the
analogous result for~$G$ a finitely generated nilpotent group of
class~two, and~${\cal N}\!il$ replaced by the category of all
metabelian nilpotent groups.

The contents of this work are part of the author's doctoral
dissertation, which was conducted under the direction of Prof.~George
M.~Bergman, at the University of California at~Berkeley. It is my very
great pleasure to express my deep gratitude and indebtedness to
Prof.~Bergman for his ever-ready advice, encouragement,
and~insights. His many suggestions and corrections have improved this
work in ways too numerous to list explicitly. Any errors that remain,
however, are my own~responsibility.

\Section{Notation and basic facts about dominions}{notation}

All groups will be written multiplicatively unless otherwise stated,
and all maps will be assumed to be group morphisms unless otherwise
specified. Given a group~$G$, the multiplicative identity of~$G$ is
denoted~$e_G$, although we will omit the subscript if it is understood
from context. Given a group~$G$, $Z(G)$ denotes the center of~$G$.

Given two groups~$A$ and~$B$, $A\wr B$ denotes the standard wreath
product of~$A$ by~$B$, that is, the semidirect product of $A^B$
by~$B$, with~$B$ acting on the index set by the regular
right~action. Elements of~$A\wr B$ will be written as $b\varphi$,
where $b\in B$ and $\varphi\in A^B$; that is, $\varphi$ is a set map
from~$B$ into~$A$.

We briefly recall some of the basic properties of dominions in
categories of~groups:

\lemma{closure}{Fix a full subcategory ${\cal C}$ of ${\cal G}roup$,
and let $G\in {\cal C}$. Then, for every subgroups $H$ and~$K$
of~$G$:}
{\parindent=30pt\it
 \item{(i)} $H\subseteq {\rm dom}_G^{\cal C}(H)$.
 \item{(ii)} If $H\subseteq K$, then ${\rm dom}_G^{\cal
C}(H)\subseteq{\rm dom}_G^{\cal C}(K)$.
 \item{(iii)} ${\rm dom}_G^{\cal C}({\rm dom}_G^{\cal C}(H)) = {\rm
dom}_G^{\cal C}(H).$\par}
{\it In particular, ${\rm dom}_G^{\cal C}(-)$ is a closure operator on the
lattice of subgroups of~$G$.}

\proof (i) follows directly from the definition of dominion. For (ii),
note that any two maps $f,g\colon G\to C$ which agree on~$K$ must also
agree on~$H$, hence on ${\rm dom}_G^{\cal C}(H)$. Finally, for (iii)
note that by the definition of dominion, the equalizer subgroups
of~$G$ which contain $H$ and those which contain ${\rm dom}_G^{\cal
C}(H)$ are exactly the~same.\endproof

\lemma{normalclosed}{If ${\cal C}$ is closed under taking quotients,
then normal subgroups are dominion closed. That is, if $N\triangleleft
G$, then ${\rm dom}_G^{\cal C}(N)=N$.}

\proof Consider the pair of maps $\pi,\zeta\colon G\to G/N$, where
$\pi$ is the canonical surjection onto the quotient, and $\zeta$ is
the zero map. Their equalizer is~$N$, so the dominion of~$N$
cannot be any~bigger.\endproof

\lemma{dirproduct}{If $G_1$, $G_2$, and $G_1\times G_2$ are groups
in~$\cal C$, and if $H_1$ is a subgroup of $G_1$, and $H_2$ a subgroup
of $G_2$, then
$${\rm dom}_{G_1\times G_2}^{\cal C}(H_1\times H_2) = {\rm
dom}_{G_1}^{\cal C}(H_1) \times {\rm dom}_{G_2}^{\cal C}(H_2).$$
That is, the dominion construction respects finite direct
products.}

\proof Identify $G_1$ with the subgroup $G_1\times \{e\}$ of
$G_1\times G_2$, and analogously with $G_2$. First we claim that
${\rm dom}_{G_1\times G_2}^{\cal C}(H_1) = {\rm dom}_{G_1}^{\cal
C}(H)\times \{e\}$. Indeed, if we compare the canonical projection
$\pi_2\colon G_1\times G_2 \to G_2$, with the zero map onto~$G_2$ we
see that the dominion of~$H_1$ must be contained in~$G_1$. Composing
any pair of maps from $G_1$ with the canonical projection $\pi_1\colon
G_1\times G_2\to G_1$ on the left, we see that the dominion cannot be
any bigger than ${\rm dom}_{G_1}^{\cal C}(H)$. And composing any maps
from $G_1\times G_2$ with the obvious immersion $i\colon G_1\to
G_1\times G_2$ on the right we see that it cannot be any
smaller. Symmetrically, 
${\rm dom}_{G_1\times G_2}^{\cal C}(H_2)= \{e\}\times{\rm
dom}_{G_2}^{\cal C}(H)$.

By \ref{closure}(ii), we have:
$${\rm dom}_{G_1}^{\cal C}(H_1)\times {\rm dom}_{G_2}^{\cal C}(H_2)
\subseteq {\rm dom}_{G_1\times G_2}^{\cal C}(H_1\times H_2).$$

For the reverse inclusion, let $(g_1,g_2)\notin {\rm dom}_{G_1\times
G_2}^{\cal C}(H_1\times H_2)$. Without loss of generality, assume that
$g_1$ is not in the dominion of~$H_1$. Therefore there exists a group
$K\in {\cal C}$ and a pair of maps 
$\phi,\psi\colon G_1\to K$
such that $\phi$ and~$\psi$ have the same restriction to~$H_1$, and
$\phi(g_1)\not=\psi(g_1)$. Considering the maps $\phi\circ \pi_1$ and
$\psi\circ\pi_1$, we see that they agree on $H_1\times H_2$, but
disagree on $(g_1,g_2)$. This proves the~lemma.\endproof

\lemma{quotients}{If ${\cal C}$ is closed under taking subgroups,
quotients, and finite direct products, then dominions respect
quotients. That is, if $G\in {\cal C}$, $H$ is a subgroup of~$G$, and
$N$ is a normal subgroup of~$G$ such that $N\subseteq H$, then
$${\rm dom}_{G/N}^{\cal C}(H/N) = \Bigl({\rm dom}_G^{\cal
C}(H)\Bigr)\Bigm/ N.$$}

\proof Trivially, ${\rm dom}_{G}^{\cal C}(H)/N\subseteq {\rm
dom}_{G/N}^{\cal C}(H/N)$, since we may compose any map $f\colon
G/N\to K$ with the canonical projection $\pi\colon G\to G/N$ to obtain
maps from $G$ to~$K$. 

For the converse inclusion, assume that $x\notin {\rm dom}_G^{\cal
C}(H)$. Therefore there exists a group $K\in {\cal C}$ and a pair of
maps $f,g\colon G\to K$ such that $f|_H=g|_H$, but $f(x)\not=g(x)$. 

Consider the induced homomorphisms $(f\times f),(f\times g)\colon
G\to K\times K$, and let $L$ be the subgroup of~$K\times K$ generated by the
images of~$G$ under these maps. Since $f$ and $g$ agree on~$H$, they
also agree on~$N$. Since $N\triangleleft G$, the common image of~$N$
under these two morphisms is normal in~$L$. This image is the subset
of the diagonal subgroup of $K\times K$ given by
$$(f\times f)(N) = \{ (f(n),f(n))\,|\, n\in N\}.$$

However, $(f\times f)(x)$ and $(f\times g)(x)$ are not in the same
coset of $(f\times f)(N)$ in~$L$. This will prove the lemma, by moding
out by $(f\times f)(N)$ to obtain maps $G/N\to L/(f\times f)(N)$ which
agree on $H/N$ but not on $xN$.

Indeed, $(f\times f)(x) (f\times g)(x)^{-1} = (e,f(x)g(x)^{-1})$, and
the second coordinate is nontrivial. Hence, this is not a diagonal
element and cannot lie in the image of~$N$. This proves the
lemma.\endproof

Dominions are closely related to group amalgams. Recall that an
{\it amalgam} of two groups $A$ and~$C$ with core~$B$ consists of
groups $A$, $C$, and~$B$, equipped with one-to-one morphisms
$\Phi_A\colon B\to A$ and $\Phi_C\colon B\to C$. We denote this
situation by writing $[A,C;B]$. We say that the amalgam is {\it weakly
embeddable} in~${\cal C}$ (where it is understood that $A$, $B$,
and~$C$ lie in~${\cal C}$) if there exist a group $M\in {\cal C}$ and
one-to-one mappings
$$\lambda_A\colon A\to M,\qquad \lambda_C\colon C\to M,\qquad
\lambda_B\colon B\to M$$
such that
$$\lambda_A\circ\Phi_A=\lambda_B,\qquad\qquad
\lambda_C\circ\Phi_C=\lambda_B.$$
We say the amalgam is {\it strongly embeddable} if, furthermore,
there is no identification between elements of $A\setminus B$ and
$C\setminus B$. Finally, by a {\it special amalgam} we mean an amalgam
$[A,C;B]$, where there exists an isomorphism $\alpha\colon A\to C$
such that \hbox{$\alpha\circ\Phi_A=\Phi_C$}. In this case, we usually
write $[A,A;B]$, with $\alpha={\rm id}_A$ being~understood.

Note that a special amalgam is always weakly embeddable. Also note
that if the amalgam $[A,A;B]$ is strongly embeddable in~$\cal C$, then
it follows that the dominion of~$B$ in~$A$ (in~${\cal C}$) is trivial,
by looking at the equalizer of the two embeddings from~$A$
into~$M$. In general, the dominion of~$B$ in~$A$ (in~${\cal C}$) is
the smallest subgroup~$D$ of~$A$, such that $B\subseteq D$, and
$[A,A;D]$ is strongly embeddable.
We refer the reader to the survey article by Higgins {\bf
[\cite{episandamalgs}]} for the~details.

Given two elements $x,y\in G$, we write $x^y=y^{-1}xy$, and
we will denote their commutator by $[x,y]=x^{-1}y^{-1}xy$. Given two
subsets $A,B$ of~$G$ (not necessarily subgroups), we denote by $[A,B]$
the subgroup of~$G$ generated by all elements $[a,b]$ with $a\in A$
and $b\in B$.  We also define inductively the left-normed commutators
of weight~\hbox{$c+1$}:
$$[x_1,\ldots,x_c,x_{c+1}] =
\bigl[[x_1,\ldots,x_c],x_{c+1}\bigr];\quad c\geq 2.$$

\defn For a group~$G$ we define the {\it lower central series} of~$G$
recursively as follows: $G_1 = G$, and $G_{c+1} = [G_c,G]$ for $c\geq
1$. We call $G_c$ the $c$-th term of the lower central series of~$G$;
$G_c$
is generated by elements of the form $[x_1,\ldots,x_c]$ for $c\geq
2$, and $x_i$ ranging over the elements of~$G$.
 
For future use we also define the first two terms of the {\it derived
series} of~$G$, $G' = [G,G]$ and $G''=[G',G']$.
Thus,~$G_2=G'=[G,G]$.

A group~$G$ is {\it nilpotent of class $c$} if and only if
$G_{c+1}=\{e\}$. A group~$G$ is {\it metabelian} if and only if
$G''=\{e\}$.

We will write ${\cal N}_c$ to denote the variety of all nilpotent
groups of class at most~$c$. We
also write ${\cal A}$ to denote the variety of all abelian groups, and
we use~${\cal A}^2$ to denote the variety of all metabelian groups
(that is, groups which are extensions of an abelian group by an
abelian~group). We write~${\cal N}\!il$ to denote the category of all
nilpotent groups, and~${\cal A}^2\cap{\cal N}\!il$ to denote the
category of all metabelian nilpotent groups. Note that the last two
categories are not varieties.

\Section{Dominions of subgroups of finitely generated groups in~${\cal
N}\!il$}{generalfg}

In this section we will prove that if~$G$ is a finitely generated
nilpotent group (of any class), and $H$ is a subgroup of~$G$, then the
dominion of~$H$ in~$G$ in the category of all nilpotent groups is
trivial. The idea for the proof is simple: first we will prove that
this is the case when~$G$ is a finite $p$-group. Using the fact that
dominions respect finite direct products, and that a finite nilpotent group
is the direct product of its Sylow subgroups, we extend the result to
the case when~$G$ is a finite nilpotent group. Finally, we use the
fact that in a finitely generated nilpotent group every subgroup is
closed in the profinite topology to extend the result to all finitely
generated nilpotent~groups.

Suppose that $G$ is a group and~$H$ is a subgroup of~$G$. If $(G_i)$
is a chief series of~$G$, the distinct terms of the series
$$(H\cap G_0, H\cap G_i,\ldots, H\cap G_n)$$
form a chief series of~$H$, which we denote by $H\cap (G_i)$. 

\thm{higmanembedd}{{\rm (G.~Higman~{\bf [\cite{higmanpgroups}]})} Let
$[A,C;B]$ be an amalgam, with $A$
and~$C$ both finite $p$-groups. The amalgam is strongly embeddable
into a finite $p$-group if and only if there exist chief series
$(A_i)$ of~$A$ and $(C_i)$ of~$C$ such that
$$B\cap (A_i) = B\cap (C_i).$$
In particular, a special amalgam of finite $p$-groups is always
strongly embeddable into a finite $p$-group.\noproof}

We deduce the following corollary:

\cor{trivforpgroups}{Let $G$ be a finite $p$-group, and let~$H$ be a
subgroup of~$G$. Let ${\cal P}$ be the category of all finite
$p$-groups. Then ${\rm dom}_G^{\cal P}(H)=H$. In particular, since
${\cal P}\subset {\cal N}\!il$, we~have
${\rm dom}_G^{{\cal N}\!il}(H) = H$\noproof}

\lemma{finiteinnilgen}{Let $G$ be a finite nilpotent group, and let~$H$
be a subgroup of~$G$. Then
${\rm dom}_G^{{\cal N}\!il}(H) = H.$}

\proof Since $G$ is finite and nilpotent, $G=\prod G_p$ where $G_p$ is
the $p$-Sylow subgroup of~$G$. Since $H$ is also nilpotent, we have
$H=\prod H_p = \prod (G_p\cap H)$. In particular, $H_p$ is a subgroup
of~$G_p$.

Since dominions respect finite direct products, we have
$${\rm dom}_G^{{\cal N}\!il}(H) = \prod \Bigl( {\rm dom}_{G_p}^{{\cal
N}\!il}\bigl(H_p\bigr)\Bigr).$$ By \ref{trivforpgroups}, each of the
dominions in the right hand side equals $H_p$, so $${\rm dom}_G^{
{\cal N}\!il}(H) = \prod H_p=H,$$ as~claimed.\endproof

Finally, to take the step from finite to finitely generated, we recall
the definition of the profinite topology of a group.

Given a group~$G$, we define the {\it profinite topology} on~$G$ to be
the coarsest topology which makes all normal subgroups of finite index
open, and makes~$G$ into a topological group (so multiplication
is a continuous map $G\times G\to G$, where $G\times G$ is given the
product topology, and the map $G\to G$ given by $g\mapsto g^{-1}$ is
also continuous).

Since the complement of a subgroup $H$ is the union of all cosets $xH$
with $x\notin H$, it follows that the normal subgroups of finite index
are both open and closed.

We say that a subgroup~$H$ of~$G$ is {\it closed} in the profinite
topology if it is a closed subset of the topological
space~$G$; equivalently, if it is the intersection of subgroups of~$G$
of finite index. Therefore, if the subgroup~$H$ is closed, and $x\in
G\setminus H$, then there exists a normal subgroup $N\triangleleft G$
such that $G/N$ is finite, and $x\notin HN$.

Recall also that a group~$G$ is said to be {\it polycyclic} iff it has
a normal series
$$\{e\}=G_n\subseteq G_{n-1}\subseteq\cdots\subseteq G_1=G$$
such that $G_{i+1}\triangleleft G_i$, and $G_i/G_{i+1}$ is cyclic. If
furthermore we may find a normal series such that $G_i\triangleleft G$
and $G_{i}/G_{i+1}$ is cyclic then~$G$ is called {\it supersolvable}.

\lemma{nilthensupersolv}{{\rm (Theorem~31.12 in~{\bf
[\cite{hneumann}]})} A finitely generated nilpotent group is
supersolvable, hence polycyclic.\noproof}

\thm{polycycthenallclosed}{{\rm (Mal'cev {\bf [\cite{maltsevtwo}]})} If~$G$
is a polycyclic group, then every subgroup of~$H$ is closed in the
profinite topology.\noproof}

In particular, every subgroup of a finitely generated nilpotent group
is~closed.

\thm{fgthentriv}{Let $G$ be a finitely generated nilpotent
group, and let $H$ be a subgroup of~$G$. Then
${\rm dom}_G^{{\cal N}\!il}(H) = H$.}

\proof Let $x\notin H$. By \ref{polycycthenallclosed} there exists a
normal subgroup $N\triangleleft G$, such that $G/N$ is finite, and
$xN\notin HN/N$. Since ${\cal N}\!il$ is closed under quotients,
subgroups, and finite direct products, it follows from \ref{quotients}
that the dominion construction in~${\cal N}\!il$ respects
quotients. Therefore,
$${\rm dom}_{G/N}^{{\cal N}\!il}(HN/N) = \Bigl({\rm dom}_G^{{\cal
N}\!il}(HN)\Bigr)\bigm/ N.\leqno(\numbeq{displayquotient})$$

By \ref{finiteinnilgen}, the dominion of~$HN/N$ in~$G/N$ is equal
to~$HN/N$. Therefore,
$$xN\notin {\rm dom}_{G/N}^{{\cal N}\!il}\bigl(HN/N\bigr).$$

By~(\ref{displayquotient}), it follows that $x\notin {\rm
dom}_G^{{\cal N}\!il}(HN)$, and hence in particular that $x$ is not
in~${\rm dom}_G^{{\cal N}\!il}(H)$. Therefore, the dominion of~$H$
in~$G$ in~${\cal N}\!il$ is equal to~$H$, as~claimed.\endproof

\Section{Dominions of subgroups of finitely generated nilpotent groups
of class~two in~${\cal A}^2\cap {\cal N}\!il$}{generalforniltwo}

In this section we will prove the analogous results to those in the
previous section, where the category of context is now~${\cal A}^2\cap
{\cal N}\!il$, and the group~$G$ is restricted to the variety of
nilpotent groups of class at most~two.

\thm{pgroups}{Let $G$ be a finite $p$-group lying in ${\cal N}_2$,
with $p$ a prime, and let~$H$ be a subgroup of~$G$. Then
${\rm dom}_G^{{\cal A}^2\cap{\cal N}\!il}(H) = H$.}

\proof Let $G$ and~$H$ be as in the statement of the lemma.
Let $N=[G,G]$ be the commutator of~$G$. First note that 
${\rm dom}_G^{{\cal A}^2\cap {\cal N}\!il}(H)\subseteq HN$.
Indeed, $HN$ is normal in~$G$, hence dominion closed by
\ref{normalclosed}. Since $H\subseteq HN$, this now follows from
\ref{closure}(ii). 

Also, note that~$N$ is
central in~$G$, since $G\in {\cal N}_2$. In particular, $N$ is abelian.
First, we define a transversal of~$N$ in~$G$ (that is, a set of coset
representatives). We claim that there is a transversal $\tau\colon
G/N\to G$ (note that~$\tau$ is only a set map, not a group morphism),
with the following properties:

{\parindent=45pt
\item{(\numbeq{proptransone})} $\tau (N) = e$.\par
\item{(\numbeq{proptransthree})}{For every $h\in H,\  y\in G$,
there exists an element $h'\in H$ such that $\tau(yh^{-1}N) =
\tau(yN)h'^{-1}$.}\par}

To construct such a transversal, consider the left action of~$H$ on
the set of cosets of~$N$, under which $h\in H$ takes $tN$ to
$tNh^{-1}=th^{-1}N$. Since~$N$ is normal, this is a well
defined action. This action partitions the set of cosets of~$N$ into
orbits. For each $H$-orbit, we first define $\tau$ to take some
arbitrary coset $tN$ in that orbit to any representative, which we now
choose once and for all, making sure to select $e$ as a representative
for~$N$. For any other coset $t'N$ in the same orbit, there exists
some element $h\in H$ such that 
$$t'\equiv \tau(tN)h^{-1}\pmod{N},$$
because this is precisely the $H$-action. Choose such an~$h$ for each
coset (the choice of~$h$ is only determined up to congruence
modulo~\hbox{$H\cap N$}), and define $\tau(t'N)=\tau(tN)h^{-1}$.

Let $\pi\colon G\to G/N$ be the canonical projection onto the
quotient. Note that $G/N$ is also abelian. For simplicity, we write
the cosets using their chosen representatives; that is, whenever we
write a coset as $tN$ it will be understood that $\tau(tN)=t$. If we
wish to represent the coset of an arbitrary element $y\in G$, where
$y$ is not the chosen representative, we will write $\pi(y)$~instead.

Since $G$ is an extension of~$N$ by $G/N$, we can embed $G$ into $N\wr
(G/N)$ by a map~$\gamma$, given by
$\gamma(g)=\pi(g)\varphi_g$,
where $\varphi_g\in N^{(G/N)}$, and for each $\pi(y)\in G/N$,
$$\varphi_g(\pi(y)) = \Bigl(\tau\bigl(y\pi(g)^{-1}\bigr) g 
\tau(\pi(y))^{-1}\Bigr)$$
(this is a theorem of Kaloujnine and Krasner {\bf [\cite{wreathext}]}).

Note that $N\wr (G/N)$ is an extension of an abelian group by an
abelian group, hence lies in~${\cal A}^2$. Since it is also a finite
$p$-group, it is nilpotent, and therefore lies in~${\cal A}^2\cap
{\cal N}\!il$.

Consider the two group morphisms $\eta,\zeta\colon N\to N/H\cap N$,
where $\eta$ is the canonical surjection, and $\zeta$ is the zero
map. The equalizer of the two maps is exactly~$H\cap N$, and
the maps induce two maps
$$\eta^*,\zeta^*\colon N\wr (G/N) \to (N/H\cap N)\wr (G/N)$$
by 
$$\eqalign{\eta^*\bigl(\pi(g)\varphi_g\bigr) &=
\pi(g)(\eta\circ\varphi_g)\cr
\zeta^*\bigl(\pi(g)\varphi_g\bigr) &= \pi(g)(\zeta\circ\varphi_g).\cr}$$

Note that $(N/H\cap N)\wr(G/N)$ is also a finite $p$-group, and metabelian.
We now consider the two maps $\eta^*\circ \gamma$ and $\zeta^*\circ\gamma$.

Let $n\in N$. By definition of $\gamma$, we have
$\gamma(n)=\varphi_n$, where $\varphi_n\colon G/N\to N$ and is given
by
$$\eqalign{\varphi_n(yN) &=
\tau\bigl(yN\pi(n)^{-1}\bigr)n\tau(yN)^{-1}\cr
&= \tau(yN)n\tau(yN)^{-1}\cr
&= yny^{-1}\cr
&= n\cr}$$
since $n\in N$, and $N$ is central. 

Therefore $\eta\circ\varphi_n(yN)=\eta(n)$, and
$\zeta\circ\varphi_n(yN)=\zeta(n)=e$. So $\eta^*\circ
\gamma(n)$ is equal to $\zeta^*\circ\gamma(n)$ if and only if
$\eta(n)=e$, that is 
if and only if $n\in H\cap N$. In particular, $(\eta^*\circ
\gamma)|_N$ and $(\zeta^*\circ\gamma)|_N$ agree exactly on~$H\cap N$.

We further claim that $\eta^*\circ\gamma$ and $\zeta^*\circ\gamma$
agree on~$H$. Indeed, let $h\in H$. Then
$\gamma(h)=\pi(h)\varphi_h$. Since $\eta^*$ and $\zeta^*$ leave the
$G/N$ component unchanged, we may concentrate on~$\varphi_h$.

Let $yN\in G/N$. By definition of $\gamma$ we have
$$\eqalign{\varphi_h(yN) & = \tau\bigl(yN\pi(h)^{-1}\bigr) h
\tau(yN)^{-1}\cr
&= \tau(yh^{-1}N)h\tau(yN)^{-1}\cr
&= yh'^{-1}hy^{-1}\cr}$$
where $yh'^{-1}= \tau (yN\pi(h)^{-1})$, with $h'\in H$. This is possible
by~(\ref{proptransthree}).

In particular, we have $h'^{-1}\equiv h^{-1} \pmod{N}$, so $h'^{-1}h\in
N$. Therefore, since $N$ is central, we have that
$\varphi_h(yN)=h'^{-1}h$, where $h'$ depends on~$y$ and~$h$, lies
in~$H$, and $h'^{-1}h\in H\cap N$.

Therefore, for every $yN\in G/N$, 
$$\eta\circ\varphi_h(yN)=\eta(h'^{-1}h)=e = \zeta(h'^{-1}h) = \zeta
\circ\varphi_h(yN).$$

Therefore the two maps agree on~$h$, and since $h$ was arbitrary, they
agree on~$H$, as~claimed. In particular,
$${\rm dom}_G^{{\cal A}^2\cap{\cal N}\!il}(H)\subseteq {\rm
Eq}(\eta^*\circ \gamma, \zeta^*\circ\gamma).$$
Therefore, ${\rm dom}_G^{{\cal A}^2\cap{\cal N}\!il}(H)\cap N = H\cap
N$.

Now consider $d\in {\rm dom}_G^{{\cal A}^2\cap {\cal
N}\!il}(H)$. Since the dominion is contained in~$HN$, there exists
$h\in H$ and $n\in N$ such that $d=hn$. In particular, $dh^{-1}=n$
lies in $N$. Since $dh^{-1}$ is also in the dominion, and is in~$N$,
it lies in $N\cap H$. In particular, $d=hn$ lies in~$H$. This proves
the required inclusion, and hence the~theorem.\endproof

The rest of the argument now proceeds as in the previous section. We
pass from finite $p$-groups to finite nilpotent groups of class~two by
decomposing the group into a direct product of its $p$-Sylow
subgroups:

\thm{finiteinmetabnil}{Let $G\in {\cal N}_2$ be a finite group
and $H$ a subgroup of~$G$. Then
${\rm dom}_G^{{\cal A}^2\cap{\cal N}\!il}(H) = H$.\noproof}

Finally, we use \ref{polycycthenallclosed} to pass from the finite case
to the finitely generated case:

\thm{fgntwothentriv}{Let $G$ be a finitely generated group lying
in~${\cal N}_2$, and let~$H$ be a subgroup of~$G$. Then
${\rm dom}_G^{{\cal A}^2\cap {\cal N}\!il}(H) =H$.\noproof}

\vskip \baselineskip
\centerline{ACKNOWLEDGEMENTS}

The author was supported in part by a fellowship from the Programa de
Formaci\'on y Superaci\'on del Personal Acad\'emico de la UNAM,
administered by the DGAPA.

%
\ifnum0<\citations{\par\bigbreak
\filbreak\centerline{\rm REFERENCES}\par\frenchspacing}\fi
%
\ifundefined{xthreeNB}\else
\item{\bf [\refer{threeNB}]}{G{.} Baumslag, B.H.~Neumann,
H.~Neumann, and P.M.~Neumann. {On varieties generated by a
finitely generated group,} {\it Math.\ Z.} {\bf 86} (1964)
\hbox{93--122}. {MR:30\#138}}\par\filbreak\fi
\ifundefined{xbergman}\else
\item{\bf [\refer{bergman}]}{George M.~Bergman. {``An Invitation to
General Algebra and Universal Constructions,''} {Berkeley
Mathematics Lecture Notes 7,} 1995.}\par\filbreak\fi
\ifundefined{xordersberg}\else
\item{\bf [\refer{ordersberg}]}{George M.~Bergman, {Ordering
coproducts of groups and semigroups,} {\it J.\ Algebra} {\bf 133} (1990)
no. 2, \hbox{313--339}. {MR:91j:06035}}\par\filbreak\fi
\ifundefined{xbirkhoff}\else
\item{\bf [\refer{birkhoff}]}{Garrett Birkhoff, {On the structure
of abstract algebras.} {\it Proc.\ Cambridge\ Philos.\ Soc.} {\bf
31} (1935), \hbox{433--454}.}\par\filbreak\fi
\ifundefined{xbrown}\else
\item{\bf [\refer{brown}]}{Kenneth S.~Brown, {``Cohomology of
Groups'' 2nd Edition,} {GTM 87\/},
Springer Verlag,~1994. {MR:96a:20072}}\par\filbreak\fi
\ifundefined{xmetab}\else
\item{\bf [\refer{metab}]}{O.N.~Golovin, {Metabelian products of
groups,}
{\it Amer.\ Math.\ Soc.\ Translations series 2}, {\bf 2} (1956),
\hbox{117--131.} {MR:17,824b}}\par\filbreak\fi
\ifundefined{xhall}\else
\item{\bf [\refer{hall}]}{M.~Hall, {``The Theory of Groups''}
Mac~Millan Company,~1959. {MR:21\#1996}}\par\filbreak\fi
\ifundefined{xphall}\else
\item{\bf [\refer{phall}]}{P.~Hall, {Verbal and marginal
subgroups,} {\it J.\ Reine\ Angew.\ Math.\/} {\bf 182} (1940)
\hbox{156--157.} {MR:2,125i}}\par\filbreak\fi
\ifundefined{xheineken}\else
\item{\bf [\refer{heineken}]}{H.~Heineken, {Engelsche Elemente der
L\"ange drei,} {\it Illinois J.\ of Math.} {\bf 5} (1961)
\hbox{681--707.} {MR:24\#A1319}}\par\filbreak\fi
\ifundefined{xherman}\else
\item{\bf [\refer{herman}]}{Krzysztof~Herman, {Some remarks on
the twelfth problem of Hanna Neumann,} {\it Publ.\ Math.\ Debrecen}
{\bf 37} (1990)  no. 1--2, \hbox{25--31.} {MR:91f:20030}}\par\filbreak\fi
\ifundefined{xherstein}\else
\item{\bf [\refer{herstein}]}{I.~N. Herstein, {``Topics in
Algebra,''} Blaisdell Publishing Co.,~1964.}\par\filbreak\fi
\ifundefined{xepisandamalgs}\else
\item{\bf [\refer{episandamalgs}]}{Peter M.~Higgins, {Epimorphisms
and amalgams,} {\it
Colloq.\ Math.} {\bf 56} no.~1 (1988) \hbox{1--17.}
{MR:89m:20083}}\par\filbreak\fi
\ifundefined{xhigmanpgroups}\else
\item{\bf [\refer{higmanpgroups}]}{Graham Higman, {Amalgams of
$p$-groups,} {\it J. of~Algebra} {\bf 1} (1964)
\hbox{301--305.} {MR:29\#4799}}\par\filbreak\fi
\ifundefined{xhigmanremarks}\else
\item{\bf [\refer{higmanremarks}]}{Graham Higman, {Some remarks
on varieties of groups,} {\it Quart.\ J.\ of Math.\ (Oxford) (2)} {\bf
10} (1959), \hbox{165--178.} {MR:22\#4756}}\par\filbreak\fi
\ifundefined{xhughes}\else
\item{\bf [\refer{hughes}]}{N.J.S.~Hughes, {The use of bilinear
mappings in the classification of groups of class~$2$,} {\it Proc.\
Amer.\ Math.\ Soc.\ } {\bf 2} (1951) \hbox{742--747.}
{MR:13,528e}}\par\filbreak\fi
\ifundefined{xisbelltwo}\else
\item{\bf [\refer{isbelltwo}]}{J.~M.~Howie and J.~R. Isbell, {
Epimorphisms and dominions II,} {\it J.\ Algebra {\bf
6}}(1967) \hbox{7--21.} {MR:35\#105b}}\par\filbreak\fi
\ifundefined{xisbellone}\else
\item{\bf [\refer{isbellone}]}{J. R. Isbell, {Epimorphisms and
dominions} {\it in} { 
``Proc.~of the Conference on Categorical Algebra, La Jolla 1965,''\/}
Lange and Springer, New
York~1966. MR:35\#105a (The statement of the
Zigzag Lemma for {\it rings} in this paper is incorrect. The correct
version is stated in~{\bf [\cite{isbellfour}]}.)}\par\filbreak\fi
\ifundefined{xisbellthree}\else
\item{\bf [\refer{isbellthree}]}{J. R. Isbell, {Epimorphisms and
dominions III,} {\it Amer.\ J.\ Math.\ }{\bf 90} (1968)
\hbox{1025--1030.} {MR:38\#5877}}\par\filbreak\fi
\ifundefined{xisbellfour}\else
\item{\bf [\refer{isbellfour}]}{J. R. Isbell, {Epimorphisms and
dominions IV,} {\it J.\ London Math.\ Soc.~(2),}
{\bf 1} (1969) \hbox{265--273.} {MR:41\#1774}}\par\filbreak\fi
\ifundefined{xjones}\else
\item{\bf [\refer{jones}]}{Gareth A.~Jones, {Varieties and simple
groups,} {\it J.\ Austral.\ Math.\ Soc.} {\bf 17} (1974)
\hbox{163--173.} {MR:49\#9081}}\par\filbreak\fi
\ifundefined{xjonsson}\else
\item{\bf [\refer{jonsson}]}{B.~J\'onsson, {Varieties of groups of
nilpotency three,} {\it Notices Amer.\ Math.\ Soc.} {\bf 13} (1966)
488.}\par\filbreak\fi
\ifundefined{xwreathext}\else
\item{\bf [\refer{wreathext}]}{L.~Kaloujnine and Marc Krasner,
{Produit complet des groupes de permutations et le probl\`eme
d'extension des groupes III,} {\it Acta Sci.\ Math.\ Szeged} {\bf 14}
(1951) \hbox{69--82}. {MR:14,242d}}\par\filbreak\fi
\ifundefined{xkhukhro}\else
\item{\bf [\refer{khukhro}]}{Evgenii I. Khukhro, {``Nilpotent Groups
and their Automorphisms,''} {de Gruyter Expositions in Mathematics}
{\bf 8}, New York 1993. {MR:94g:20046}}\par\filbreak\fi
\ifundefined{xkleimanbig}\else
\item{\bf [\refer{kleimanbig}]}{Yu.~G. Kle\u{\i}man, {On
identities in groups,} {\it Trans.\ Moscow Math.\ Soc.\ } 1983,
Issue 2, \hbox{63--110}. {MR:84e:20040}}\par\filbreak\fi
\ifundefined{xthirtynine}\else
\item{\bf [\refer{thirtynine}]}{L. G. Kov\'acs, {The thirty-nine
varieties,} {\it Math.\ Scientist} {\bf 4} (1979)
\hbox{113--128.} {MR:81m:20037}}\par\filbreak\fi
\ifundefined{xlamssix}\else
\item{\bf [\refer{lamssix}]}{T.Y. Lam and David B. Leep, {
Combinatorial structure on the automorphism group of~$S_6$,} {\it
Expo. Math.} {\bf 11} (1993) \hbox{289--308.}
{MR:94i:20006}}\par\filbreak\fi
\ifundefined{xlevione}\else
\item{\bf [\refer{levione}]}{F.~W. Levi, {Groups on which the
commutator relation 
satisfies certain algebraic conditions,} {\it J.\ Indian Math.\ Soc.\ New
Series} {\bf 6}(1942), \hbox{87--97.} {MR:4,133i}}\par\filbreak\fi
\ifundefined{xgermanlevi}\else
\item{\bf [\refer{germanlevi}]}{F.~W. Levi and B. L. van der Waerden,
{\"Uber eine 
besondere Klasse von Gruppen,} {\it Abhandl.\ Math.\ Sem.\ Univ.\ Hamburg}
{\bf 9}(1932), \hbox{154--158.}}\par\filbreak\fi
\ifundefined{xlichtman}\else
\item{\bf [\refer{lichtman}]}{A. L. Lichtman, {Necessary and
sufficient conditions for the residual nilpotence of free products of
groups,} {\it J. Pure and Applied Algebra} {\bf 12} no. 1 (1978),
\hbox{49--64.} {MR:58\#5938}}\par\filbreak\fi
\ifundefined{xmaxofan}\else
\item{\bf [\refer{maxofan}]}{Martin W. Liebeck, Cheryl E. Praeger, 
and Jan Saxl, {A classification of the maximal subgroups of the
finite alternating and symmetric groups,} {\it J. of Algebra} {\bf
111}(1987), \hbox{365--383.} {MR:89b:20008}}\par\filbreak\fi
\ifundefined{xepisingroups}\else
\item{\bf [\refer{episingroups}]}{C.E. Linderholm, {A group
epimorphism is surjective,} {\it Amer.\ Math.\ Monthly\ }77
\hbox{176--177.}}\par\filbreak\fi
\ifundefined{xmckay}\else
\item{\bf [\refer{mckay}]}{Susan McKay, {Surjective epimorphisms
in classes
of groups,} {\it Quart.\ J.\ Math.\ Oxford (2),\/} {\bf 20} (1969),
\hbox{87--90.} {MR:39\#1558}}\par\filbreak\fi
\ifundefined{xmachenry}\else
\item{\bf [\refer{machenry}]}{T. MacHenry, {The tensor product and
the 2nd nilpotent product of groups,} {\it Math. Z.\/} {\bf 73}
(1960), \hbox{134--145.} {MR:22\#11027a}}\par\filbreak\fi

\ifundefined{xmaclane}\else
\item{\bf [\refer{maclane}]}{Saunders Mac Lane, {``Categories for
the Working Mathematician,''} {GTM 5},
Springer Verlag (1971). {MR:50\#7275}}\par\filbreak\fi
\ifundefined{xbilinearprelim}\else
\item{\bf [\refer{bilinearprelim}]}{Arturo Magidin, {Bilinear maps
and central extensions of abelian groups,} {\it Submitted.}}\par\filbreak\fi
\ifundefined{xprodvarprelim}\else
\item{\bf [\refer{prodvarprelim}]}{Arturo Magidin, {Dominions in decomposable
varieties of groups,} {\it Submitted.}}\par\filbreak\fi
\ifundefined{xmythesis}\else
\item{\bf [\refer{mythesis}]}{Arturo Magidin, {``Dominions in
Varieties of Groups,''} Doctoral dissertation, University of
California at Berkeley, May 1998.}\par\filbreak\fi
\ifundefined{xnildomsprelim}\else
\item{\bf [\refer{nildomsprelim}]}{Arturo Magidin, {Dominions in
varieties of nilpotent groups,} {\it Comm.\ Alg.} to appear.}\par\filbreak\fi
\ifundefined{xsimpleprelim}\else
\item{\bf [\refer{simpleprelim}]}{Arturo Magidin, {Dominions in
varieties generated by simple groups,} {\it Submitted.}}\par\filbreak\fi
\ifundefined{xdomsmetabprelim}\else
\item{\bf [\refer{domsmetabprelim}]}{Arturo Magidin, {Dominions
in the variety of metabelian groups,}
{\it Submitted.}}\par\filbreak\fi
\ifundefined{xfgnilprelim}\else
\item{\bf [\refer{fgnilprelim}]}{Arturo Magidin, {Dominions of
finitely generated nilpotent groups,} {\it
Comm.\ Algebra,} to appear}\par\filbreak\fi 
\ifundefined{xwordsprelim}\else
\item{\bf [\refer{wordsprelim}]}{Arturo Magidin, {
Words and dominions,} {\it Submitted.}}\par\filbreak\fi
\ifundefined{xabsclosed}\else
\item{\bf [\refer{absclosed}]}{Arturo Magidin, {Absolutely closed
nil-2 groups,} {\it Submitted.}}\par\filbreak\fi
\ifundefined{xmagnus}\else
\item{\bf [\refer{magnus}]}{Wilhelm Magnus, Abraham Karras, and
Donald Solitar, {``Combinatorial Group Theory,''} 2nd Edition; Dover
Publications, Inc.~1976. {MR:53\#10423}}\par\filbreak\fi
\ifundefined{xamalgtwo}\else
\item{\bf [\refer{amalgtwo}]}{Berthold J. Maier, {Amalgame
nilpotenter Gruppen
der Klasse zwei II,} {\it Publ.\ Math.\ Debrecen} {\bf 33}(1986),
\hbox{43--52.} {MR:87k:20050}}\par\filbreak\fi
\ifundefined{xnilexpp}\else
\item{\bf [\refer{nilexpp}]}{Berthold J. Maier, {On nilpotent
groups of exponent $p$,} {\it J.\ Algebra} {\bf 127} (1989)
\hbox{279--289.} {MR:91b:20046}}\par\filbreak\fi
\ifundefined{xmaltsev}\else
\item{\bf [\refer{maltsev}]}{A. I. Maltsev, {Generalized
nilpotent algebras and their associated groups} (Russian), {\it
Mat.\ Sbornik N.S.} {\bf 25(67)} (1949) \hbox{347--366.} ({\it
Amer.\ Math.\ Soc.\ Translations Series 2} {\bf 69} 1968,
\hbox{1--21.}) {MR:11,323b}}\par\filbreak\fi
\ifundefined{xmaltsevtwo}\else
\item{\bf [\refer{maltsevtwo}]}{A. I. Maltsev, {Homomorphisms onto
finite groups} (Russian), {\it Ivanov. gosudarst. ped. Inst., u\v
cenye zap., fiz-mat. Nauk} {\bf 18} (1958)
\hbox{49--60.}}\par\filbreak\fi
\ifundefined{xmorandual}\else
\item{\bf [\refer{morandual}]}{S. Moran, {Duals of a verbal
subgroup,} {\it J.\ London Math.\ Soc.} {\bf 33} (1958)
\hbox{220--236.} {MR:20\#3909}}\par\filbreak\fi
\ifundefined{xhneumann}\else
\item{\bf [\refer{hneumann}]}{Hanna Neumann, {``Varieties of
Groups,''} {Ergebnisse der Mathematik und ihrer Grenz\-ge\-biete,\/}
New series, Vol.~37, Springer Verlag~1967. {MR:35\#6734}}\par\filbreak\fi
\ifundefined{xneumannwreath}\else
\item{\bf [\refer{neumannwreath}]}{Peter M.~Neumann, {On the
structure of standard wreath products of groups,} {\it Math.\
Zeitschr.\ }{\bf 84} (1964) \hbox{343--373.} {MR:32\#5719}}\par\filbreak\fi
\ifundefined{xpneumann}\else
\item{\bf [\refer{pneumann}]}{Peter M.~Neumann, {Splitting groups
and projectives
in varieties of groups,} {\it Quart.\ J.\ Math.\ Oxford} (2), {\bf
18} (1967),
\hbox{325--332.} {MR:36\#3859}}\par\filbreak\fi
\ifundefined{xoates}\else
\item{\bf [\refer{oates}]}{Sheila Oates, {Identical Relations in
Groups,} {\it J.\ London Math.\ Soc.} {\bf 38} (1963),
\hbox{71--78.} {MR:26\#5043}}\par\filbreak\fi
\ifundefined{xolsanskii}\else
\item{\bf [\refer{olsanskii}]}{A. Ju.~Ol'\v{s}anski\v{\i}, {On the
problem of a finite basis of identities in groups,} {\it
Izv.\ Akad.\ Nauk.\ SSSR} {\bf 4} (1970) no. 2
\hbox{381--389.}}\par\filbreak\fi
\ifundefined{xremak}\else
\item{\bf [\refer{remak}]}{R. Remak, {\"Uber minimale invariante
Untergruppen in der Theorie der end\-lichen Gruppen,} {\it
J.\ reine.\ angew.\ Math.} {\bf 162} (1930),
\hbox{1--16.}}\par\filbreak\fi
\ifundefined{xclassifthree}\else
\item{\bf [\refer{classifthree}]}{V. N. Remeslennikov, {Two
remarks on 3-step nilpotent groups} (Russian), {\it Algebra i Logika
Sem.} (1965) no.~2 \hbox{59--65.} {MR:31\#4838}}\par\filbreak\fi
\ifundefined{xrotman}\else
\item{\bf [\refer{rotman}]}{J.J. Rotman, {``Introduction to the Theory of
Groups 4th edition,''} {GTM~119},
Springer Verlag,~1994. {MR:95m:20001}}\par\filbreak\fi
\ifundefined{xsaracino}\else
\item{\bf [\refer{saracino}]}{D. Saracino, {Amalgamation bases for
nil-$2$ groups,} {\it Alg.\ Universalis\/} {\bf 16} (1983),
\hbox{47--62.} {MR:84i:20035}}\par\filbreak\fi
\ifundefined{xscott}\else
\item{\bf [\refer{scott}]}{W. R. Scott, {``Group Theory,''} Prentice
Hall,~1964. {MR:29\#4785}}\par\filbreak\fi
\ifundefined{xsmelkin}\else
\item{\bf [\refer{smelkin}]}{A. L. \v{S}mel'kin, {Wreath products and
varieties of groups} (Russian), {\it Dokl.\ Akad.\ Nauk S.S.S.R.\/} {\bf
157} (1964), \hbox{1063--1065} Transl.: {\it Soviet Math.\ Dokl.\ } {\bf
5} (1964), \hbox{1099--1011}. {MR:33\#1352}}\par\filbreak\fi
\ifundefined{xstruikone}\else
\item{\bf [\refer{struikone}]}{Ruth Rebekka Struik, {On nilpotent
products of cyclic groups,} {\it Canadian J.\ Math.}
{\bf 12} (1960)
\hbox{447--462}. {MR:22\#11028}}\par\filbreak\fi
\ifundefined{xstruiktwo}\else
\item{\bf [\refer{struiktwo}]}{Ruth Rebekka Struik, {On nilpotent
products of cyclic groups II,} {\it Canadian J.\
Math.\/} {\bf 13} (1961) \hbox{557--568.}
{MR:26\#2486}}\par\filbreak\fi
\ifundefined{xvlee}\else
\item{\bf [\refer{vlee}]}{M. R. Vaughan-Lee, {Uncountably many
varieties of groups,} {\it Bull.\ London Math.\ Soc.} {\bf 2} (1970)
\hbox{280--286.} {MR:43\#2054}}\par\filbreak\fi
\ifundefined{xweibel}\else
\item{\bf [\refer{weibel}]}{Charles Weibel, {``Introduction to
Homological Algebra,''} Cambridge University
Press~1994. {MR:95f:18001}}\par\filbreak\fi 
\ifundefined{xweigelone}\else
\item{\bf [\refer{weigelone}]}{T. S. Weigel, {Residual properties
of free groups,} {\it J.\ Algebra} {\bf 160} (1993)
\hbox{14--41.} {MR:94f:20058a}}\par\filbreak\fi
\ifundefined{xweigeltwo}\else
\item{\bf [\refer{weigeltwo}]}{T. S. Weigel, {Residual properties
of free groups II,} {\it Comm.\ in Algebra} {\bf 20}(5) (1992)
\hbox{1395--1425.} {MR:94f:20058b}}\par\filbreak\fi
\ifundefined{xweigelthree}\else 
\item{\bf [\refer{weigelthree}]}{T. S. Weigel, {Residual Properties
of free groups III,} {\it Israel J.\ Math.\ } {\bf 77} (1992)
\hbox{65--81.} {MR:94f:20058c}}\par\filbreak\fi
\ifundefined{xzstwo}\else
\item{\bf [\refer{zstwo}]}{Oscar Zariski and Pierre Samuel,
{``Commutative Algebra, Volume
II,''} Springer-Verlag~1976. {MR:52\#10706}}\par\filbreak\fi
\ifnum0<\citations\nonfrenchspacing\fi

\bigskip

\vfill
\eject
\immediate\closeout\aux
\end